\documentclass{amsart}

\newcommand{\EQREFER}{\ref}
\newcommand{\EQXREF}{\label}
\newcommand{\REFER}{\ref}

\newcommand{\C}{{\mathbb C}}

\newcommand{\CC}{{\overline \C}}

\newcommand{\A}{{\mathcal A}}
\newcommand{\downto}{{\searrow}}
\newcommand{\F}{{\mathcal F}}
\newcommand{\E}{{\mathcal E}}
\newcommand{\nbhd}{{neighborhood }}
\newcommand{\tri}{{\bigtriangleup}}

\newtheorem{corollary}{Corollary}
\newtheorem{lemma}{Lemma}
\newtheorem{property}{Property}
\newtheorem{theorem}{Theorem}

\newtheorem{othertheorem}{Theorem}

\theoremstyle{definition}
\newtheorem{definition}{Definition}
\theoremstyle{remark}
\newtheorem{remark}{Remark}

\numberwithin{equation}{section}

\begin{document}

\title{Completely Invariant Julia sets of Polynomial Semigroups}
\author{Rich Stankewitz}
\thanks{Research supported by a Department of Education 
GAANN fellowship and by the Research
Board of the University of Illinois at Urbana-Champaign.}
\address{Department of Mathematics, 
         University of Illinois, 
         Urbana, Illinois 61801}
\curraddr{Department of Mathematics, Texas A\&M University,
  College Station, TX 77843}
\email{richs@math.tamu.edu}

\begin{abstract}
Let $G$ be a semigroup of rational functions of degree at least two,
under composition of functions.  Suppose that $G$ contains two
polynomials with non-equal Julia sets.  We prove that the smallest
closed subset of the Riemann sphere which contains at least three
points and is completely invariant under each element of $G$, is the
sphere itself.
\end{abstract}

\maketitle

\section{Introduction}
This paper addresses the issue of how one can extend the
definition of a Julia set of a rational function of a
complex variable of degree at least two 
to the more general setting of a Julia set of a
rational semigroup.

It is possible to define the Julia set of a single rational function 
in two
different, yet equivalent ways.  The first definition of a Julia set
is given as the complement of the set of normality of the iterates of
the given function.  The second definition of a Julia set is given as
the smallest closed completely invariant set which contains three or
more points.  Each of these definitions can be naturally extended to
the setting of an arbitrary rational semigroup $G$, but 
the extensions are not equivalent.  

This paper will show that the extension of the first definition,
 denoted $J(G)$, 
is better for the purpose of achieving a situation where it 
is meaningful to study dynamics on the components of its complement, 
the Fatou set.  (If one is studying dynamics from the point of
view that complete invariance is required, then, of course, the
 extension of the second
definition, denoted $E(G)$, is better.)  We prove that for a 
semigroup generated by two polynomials, of degree greater than or equal
to two, with non-equal Julia sets, we have $E(G)=\CC$, while 
$J(G)$ is known to be a
compact subset of the complex plane $\C$.
(If the Julia sets of the two
generators are equal, then both $J(G)$ and $E(G)$
are equal to this common Julia set.)  

\section{Definitions and basic facts}

In what follows all notions of convergence will be with respect to the 
spherical metric on the Riemann sphere $\CC.$

 

A rational semigroup $G$ is a semigroup 
of rational functions of degree greater than or equal to two
defined on the Riemann sphere $\CC$ with the semigroup operation being
functional composition.  When a semigroup $G$ is generated by the functions
$\{f_1, f_2, \dots, f_n, \dots\}$,
we write this as
\begin{equation}
G=\langle f_1, f_2,\dots, f_n, \ldots \rangle.\notag
\end{equation}

In~\cite{HM1}, p.~360 the definitions of the set of
normality, often called the Fatou set, 
and the Julia set of a rational semigroup are as follows:

\begin{definition} \label{N(G),J(G)}
For a rational semigroup $G$ we
define the set of normality of $G$, $N(G)$, by 
$$N(G)=\{z \in \CC:\exists \text{ a \nbhd of } z \text{ on
which } G \text{ is a normal family} \}$$
and define the Julia set of $G$, $J(G)$, by
$$J(G)=\CC \setminus N(G).$$
\end{definition}

Clearly from these definitions we see that $N(G)$ is an open set and
therefore its complement $J(G)$ is a compact set.
These definitions generalize the case of iteration of a single rational
function and we write $N(\langle h \rangle )=N_h$ and
$J(\langle h \rangle )=J_h$.  

Note that $J(G)$ contains the Julia set of each element of $G$.

\begin{definition} If $h$ is a map of a set $Y$ into
itself, a subset $X$ of $Y$ is:
\begin{align} 
  i)& \,forward\,\,invariant \text{ under } h \text{ if } h(X) \subset X;\notag\\
 ii)& \,backward\,\,invariant \text{ under } h \text{ if }  
h^{-1}(X) \subset X;\notag\\
iii)& \,completely\,\,invariant \text{ under } h \text{ if } h(X) \subset
  X \text{ and } h^{-1}(X) \subset X.\notag
\end{align}

\end{definition}

It is well known that the set of normality of $h$ and the Julia set of
$h$ are
completely invariant under $h$ (see
~\cite{Be}, p.~54), i.e., 
\begin{equation}\label{cinv}
h(N_h)=N_h=h^{-1}(N_h) \text{ and }
h(J_h)=J_h=h^{-1}(J_h).
\end{equation}

In fact we have the following result.

\begin{property}\label{Jh2}
The set $J_h$ is the smallest closed completely invariant (under $h$)
set which contains three or more points (see ~\cite{Be}, p.~67).  
\end{property}
In fact, this may be chosen
as an alternate definition to the definition of $J_h$ 
given in Definition~\ref{N(G),J(G)}.

From Definition~\ref{N(G),J(G)}, we get that $N(G)$ is 
forward invariant under
each element of $G$ and $J(G)$ is backward invariant under each
element of $G$ (see ~\cite{HM1}, p.~360). 
The sets $N(G)$ and $J(G)$ are, however, not necessarily 
completely invariant under the elements of $G$.  
This is in contrast to the case of single function dynamics as noted
in (\ref{cinv}). 
The question then arises, what if we required 
the Julia set of the semigroup $G$ to be completely 
invariant under each element of $G$?  
We consider in this paper the consequences
of such an extension, given in the following
definition. 

\begin{definition}\label{Edef}
 For a rational semigroup $G$ we define
$$E=E(G)=\bigcap\{S:S \text{ is closed, completely invariant under
each }g \in G, \#(S)\geq3 \}$$
where $\#(S)$ denotes the cardinality of $S$.
\end{definition}

We note that $E(G)$ exists, is closed, 
is completely invariant under each element
of $G$ and contains the Julia set of each
element of $G$ by Property~\REFER{Jh2}.

\vskip.1truein
We now compare the sets $E(G)$ and $J(G)$. 
\vskip.1truein

{\bf{Example 1}}
 Suppose that $G=\langle f,g \rangle$
and $J_f=J_g$.
Then $E=J_f=J_g$ since $J_f$ is completely invariant under
$f$ and $J_g$ is completely invariant under $g$.  It is easily verified that
if $J_f=J_g$, then $J(G)=J_f=J_g$, also.  

\vskip.1truein

The following
example shows that we may have 
$J(G) \ne E(G)$.

\vskip.1truein

{\bf{Example 2}}
Let $a \in \CC, |a|>1$ and $G=\langle z^2,z^2/a \rangle .$  One can easily show
that $J(G)=\{z:1 \leq |z| \leq |a| \}$ (see ~\cite{HM1}, p.360) while
$E=\CC$.  Note that $J_{z^2}=\{z:|z|=1\}$ and $J_{z^2/a}=\{z:|z|=|a|\}.$

\vskip.1truein

The main results of this paper will show that, in some sense, Examples
1 and 2 illustrate the only two possibilities for polynomial
semigroups.  
We prove the following theorems.

\begin{theorem}\label{main}
For polynomials $f$ and $g$ of degree
greater than or equal to two,
$J_f \neq J_g$ implies $E(G)=\CC$ where $G=\langle f,g \rangle.$
\end{theorem}

The following theorem follows immediately.

\begin{theorem}\label{main2}
For a rational semigroup $G'$ which contains two polynomials 
$f$ and $g$ of degree
greater than or equal to two, $J_f \neq J_g$ implies $E(G')=\CC$.
\end{theorem}

\section{Auxiliary results}

Suppose that 
$G$ is a semigroup generated by two polynomials, $f$ and $g$, 
of degree $k$ and $l$ where $k,l \geq 2$.
Note that $\infty \notin J(G)$, since any small \nbhd of
$\infty$ which is forward invariant under both $f$ and $g$ is
necessarily in $N(G)$.  This $G$
will remain fixed for the rest of the paper.  Recall that $E(G) \supset J_f
\cup J_g$.

We first begin with a lemma that will give a convenient description of
$E$.  Our next result shows how $E$ is ``built up'' from $J_f$
and $J_g$.

For a collection of sets $\A$, and a function $h$, we denote new
collections of sets by 
$h(\A)=\{h(A):A \in \A\}$ and $h^{-1}(\A)=\{h^{-1}(A):A \in \A\}$.

Let us define the following countable collections of sets:
\begin{align}
&\E_0=\{ J_f, J_g\}, \notag\\
&\E_{n+1}= f^{-1}(\E_n)  \cup  f(\E_n)  \cup
 g^{-1}(\E_n)  \cup  g(\E_n) , \notag\\
\text{ and }&\E = \bigcup_{n=0}^\infty \E_n.\notag
\end{align}

Since $E$ is completely invariant under $f$ and $g$ and contains both
$J_f$ and $J_g$, we
have  
$E \supset \bigcup_{A \in \E} A$.  Since $E$ is also closed, we have

\begin{equation}\label{E>union}
E \supset \overline{\bigcup_{A \in \E} A}, 
\end{equation}
where $\overline{B}$ denotes the closure of the set $B$.
The following
lemma shows that these two sets are actually equal.

\begin{lemma}\label{E}
We have $E=\overline{\bigcup_{A \in \E} A}.$
\end{lemma}

\begin{proof}
Because of (\ref{E>union}) we need only show that
$E \subset \overline{\bigcup_{A \in \E} A}$.
Since the set on the right is closed and contains both $J_f$ and
$J_g$, it remains only to show that it is also completely invariant
under both $f$ and $g$.  This can easily be shown by using the fact that
$f$ and $g$ are continuous open maps.
\end{proof}





\begin{corollary}\label{Eperf}
The set $E$ has no isolated points; i.e., $E$ is perfect.
\end{corollary}

\begin{proof}
Since $J_f$ and $J_g$ are perfect
(see ~\cite{Be}, p.~68) and
backward and forward images of perfect sets under rational maps are
perfect, we see that each set in $\E$ is perfect by a routine
inductive argument.  The corollary then follows since the closure of a
union of perfect sets is perfect.
\end{proof}

\begin{lemma}\label{Eint}
 If $E$ has nonempty interior, then $E=\CC.$
\end{lemma}

\begin{proof}
Suppose that $E^o \neq \emptyset,$ where $B^o$ denotes the interior of
the set $B$.
It will be shown that this implies the existence of an open set $U
\subset E$ such that $U$ intersects $J_f $ or $J_g$.  Supposing
that $U$ intersects $J_g$, we may observe by
the expanding property of Julia sets (see ~\cite{Be}, p.~69) that we have
$\overline{\bigcup_{n=1}^\infty g^n(U)} = \CC.$  By the forward
invariance of the closed set $E$ under the map $g$ we see that this 
implies that $E=\CC$.

We will use the following elementary fact:  

For any sets $B$ and $C$ and
any function $h$ we have
\begin{equation}\label{fact}
 B \cap h(C) \neq \emptyset \text{ if and only if }
h^{-1}(B) \cap C \neq \emptyset.
\end{equation}

Since $E^o \neq \emptyset$, there exists an open disc  $\tri
\subset E$.  By Lemma~\REFER{E} we see then that there exists a set $A$
in $\E_n$, say,  such that $\tri \cap A \neq \emptyset.$  Since $A \in
\E_n$, we see that it can be expressed as $ A=h_n \cdots h_1(J_g)$,
for example,
where each $h_j \in \{f,f^{-1},g,g^{-1} \}$.  Considering each $h_j$
as a map on subsets of $\CC$, as opposed to a map on points of $\CC$,
we can define the ``inverse'' maps $h_j^*$
accordingly, i.e., $h_1=f$ implies $h_1^*=f^{-1}$ and $h_2=g^{-1}$
implies $h_2^*=g$.  The $h_j^*$ are not true inverses since
$f^{-1}(f(X))$ may properly contain $X$.

The fact (\ref{fact}) does imply, however, that 
\begin{align}
 A \cap \tri \neq \emptyset
 \implies & 
 h_n \cdots h_1(J_g)\cap \tri \neq \emptyset \notag\\
 \implies & 
h_{n-1} \cdots h_1(J_g)\cap h_n^*(\tri) \neq \emptyset \notag\\
& \vdots \notag\\
 \implies & h_1(J_g) \cap h_2^* \cdots h_n^*(\tri) \neq \emptyset \notag\\ 
\implies & J_g \cap  h_1^* \cdots h_n^*(\tri) \neq \emptyset. \notag
\end{align}

Since each $h_j^*$ maps open sets to open sets (as $f,g,f^{-1}, g^{-1}$
do) we see that $U=h_1^* \cdots h_n^*(\tri)$ is open.  By the complete
invariance of $E$ under both $f$ and $g$ we see that each $h_j^*$ takes
subsets of $E$ to subsets of $E$.  Hence $U \subset E$ and the proof
of Lemma~\ref{Eint} is complete.
\end{proof}

Similar to the description of the set $E$ given in Lemma~\REFER{E}
is the following description of the
Julia set $J(G)$ of the semi-group $G=\langle f,g \rangle$. 
Consider the countable collection of sets  
\begin{align}
&\F_0=\{ J_f, J_g\}, \notag\\
&\F_{n+1}=  f^{-1}(\F_n) \cup  g^{-1}(\F_n), \notag\\
\text{ and }&\F = \bigcup_{n=0}^\infty \F_n .\notag
\end{align}

Since $J(G)$ is backward invariant under $f$ and $g$, closed, and
contains both $J_f$ and $J_g$ (see ~\cite{HM1}, p.~360),  we
have $J(G) \supset \overline{\bigcup_{A \in \F} A}$.  

\begin{lemma}\label{J(G)=} 
We have $J(G)=\overline{\bigcup_{A \in \F} A}.$
\end{lemma}

\begin{proof}
Since the set on the right is closed, backward invariant for both $f$
and $g$ (as follows as in the proof of Lemma~\REFER{E}) and clearly contains
more that three points, it must contain $J(G)$ as the complement is
then in the set of normality of $G$.
\end{proof}

From Lemma~\REFER{J(G)=}, noting that $\F \subset \E$, and Lemma~\REFER{E}
we get the following result.

\begin{corollary}\label{J(G)<E}
The Julia set of the semigroup $G$ is contained
in $E$, i.e., $J(G) \subset E$.
\end{corollary}

Hence, combining Lemma~\REFER{Eint} and Corollary~\REFER{J(G)<E}, 
we get the following corollary.

\begin{corollary} 
If $J(G)$ has nonempty interior, then $E=\CC$.
\end{corollary}

\section{Proof of the main result}

In this section we will prove the main result, Theorem~\REFER{main2},
but first we establish the necessary lemmas.

\begin{lemma}\label{inftyinE}
If $f$ and $g$ are polynomials of degree
greater than or equal to two and $J_f \neq J_g$,
then $\infty \in E$.
\end{lemma}

\begin{proof}
Denoting the unbounded components of the respective Fatou sets of $f$
and $g$ by $F_\infty$ and $G_\infty$, we recall (see ~\cite{Be}, p.~54 and
p.~82) that
$J_f=\partial F_\infty$ and $J_g=\partial G_\infty$.  

Since $F_\infty$ and $G_\infty$ are domains with nonempty intersection
and $\partial F_\infty \neq \partial G_\infty$, we have
$J_f \cap G_\infty \neq \emptyset$ or 
$J_g \cap F_\infty \neq \emptyset.$ 

Hence we may select $z \in J_g \cap F_\infty$, say. 
Denoting the $n$th iterate of $f$ by $f^n$, we see 
that $f^n(z) \to \infty,$ and by
the forward invariance under the map $f$ of the set $E$ we get that
each $f^n(z) \in E$.  Since $E$ is closed we see then that
$\infty \in E$.  
\end{proof}

\begin{remark}\label{JfneqJg}
Since it will be necessary later, we make special note of 
the fact used in the above proof that 
$J_f \neq J_g$
implies $J_f \cap G_\infty \neq \emptyset$ or 
$J_g \cap F_\infty \neq \emptyset.$
\end{remark}

\begin{remark}\label{inftynotisol}
Note that the proof above shows also that $\infty$ is not an isolated
point of $E$ when $J_f \neq J_g$.  This, of course, also follows from
Corollary~\REFER{Eperf} and Lemma~\ref{inftyinE}.  
\end{remark}

The disc centered at the point $z$ with radius $r$ will be denoted
 $\tri(z,r)$.

\begin{lemma}\label{circles}
Suppose that $\tri(0,r^*)=A \cup B$ where $A$ is open, $A$ and $B$ are
disjoint, and both $A$ and $B$ are 
nonempty.  If both $A$ and $B$ are 
completely invariant under the map $L(z) = z^j$ defined
on $\tri(0,r^*)$ where $0<r^*<1$ and $j \geq 2$, then the 
set $A$ is a union of open annuli centered
at the origin and hence $B$ is a union of circles centered at the
origin.  Furthermore, each of $A$ and $B$ contains a sequence of circles
tending to zero.
\end{lemma}

\begin{proof}
Let $z_0=re^{i\theta} \in A$.  Since $A$ is open we may choose
$\delta > 0$
such that the arc $ \alpha_{z_0}=\{re^{i\omega}:|\theta-\omega| 
\leq {\frac{\delta} {2}} \} \subset A.$

Fix a positive integer $n$ such that $j^n \delta>2\pi$.  Since
$L^n(z)=z^{j^n}$ we get
\begin{equation*}
L^n(\alpha_{z_0}) = C(0,r^{j^n})
\end{equation*}
where $C(z,r)=\{\zeta:|\zeta -z|=r\}.$

By the forward invariance of $A$ under
$L$, we see that $C(0,r^{j^n}) \subset A$.  But now
by the backward invariance of $A$, 
we get
\begin{equation*}
C(0,r)=L^{-n}(C(0,r^{j^n})) \subset A.
\end{equation*}

Thus for any
$re^{i\theta} \in A$, we have $C(0,r) \subset A.$  Hence  $A$, being
open, must be a union of open annuli
centered at the origin and therefore $B$, being the complement of
$A$ in $\tri(0,r^*)$, must be a union of circles
centered at the origin.   

We also note that if $C(0,r) \subset A$, then $C(0,r^{j^n}) \subset A$ is a
sequence of circles tending to zero.  Similarly we obtain a sequence
of circles in $B$ tending to zero.
\end{proof}

\begin{lemma}\label{L=azj}  
Let $L:\tri(0,r^*) \to \tri(0,r^*)$, where $0<r^*<1$, be an analytic
function such that $L(0)=0$.
Let $B$ be a set with empty interior which is a 
union of circles centered at the origin and which contains a
sequence of circles tending to zero.  If $B$ is forward invariant under
the map $L$, then $L$ is of the form
$$L(z)=az^j$$ 
for some non-zero complex number $a$ and some positive integer $j$. 
\end{lemma}

\begin{proof}
Since $L(0)=0$, we have, near $z=0$, 
\begin{align} 
L(z)
&= az^j + a_1z^{j+1} + \cdots \notag\\ 
&= az^j[1+ {\frac {a_1} {a}} z+ \cdots] \notag
\end{align}
for some non-zero complex number $a$ and some positive integer $j$. 

Let $h(z)=L(z)/az^j$ and note that $h(z)$ is analytic and
tends to  1 as $z$ tends to
$0.$  We shall prove that $h(z) \equiv 1$ and the lemma then follows.

Let $C_n=C(0,r_n)$ be sequence of circles contained in $B$ with $r_n \to
0.$   We claim that each $L(C_n)$ is contained in another
circle centered at the origin of, say, radius $r_n'.$  If not, then
the connected set
$L(C_n)$ would contain points of all moduli between, say, $r'$ and
$r''$.  This, however, would imply that $B$
would contain the annulus between the circles $C(0,r')$ and
$C(0,r'')$. Thus we have $L(C_n) \subset C(0,r_n')$.

So we see then that $h(C_n) \subset C(0,r_n'/|a|r_n^j).$
  
But for large $n$ we see that if $h$ were non-constant, then $h(C_n)$ would
be a path which stays near $h(0)=1$ and winds around $h(0)=1$.  
Since $h(C_n)$ is contained in a circle
centered at the origin, this cannot happen.  We thus conclude that $h$
is constant.
\end{proof}

\begin{lemma}\label{|a|<1}
If $B \subset \tri(0,r^*)$ for $0<r^*<1$ 
is a nonempty relatively closed set 
which is completely
invariant under the maps $H:z \mapsto z^j$ and $K:z \mapsto az^m$
defined on $\tri(0,r^*)$ where $a$ is a nonzero complex number and
$j,m$ are integers with $j,m \geq 2$, 
then $B=\tri(0,r^*)$ or $|a|=1.$
\end{lemma}

\begin{proof}
We may assume that $|a| \leq 1$ by the following reasoning.
Suppose that $|a| \geq 1.$  Let $b$ be a complex number such that $b^{m-1}=a$
and define $\psi(z)=bz.$  Since $\psi \circ H \circ \psi^{-1}(z) =
z^j/b^{j-1}$ and $\psi \circ K \circ \psi^{-1}(z) = z^m$,
we see that the lemma would then imply that $\psi(B)=\tri(0,|b|r^*)$ 
or $|b|=1$.  Since we know that $\psi(B)=\tri(0,|b|r^*)$ 
exactly when $B=\tri(0,r^*)$, and $|b|=1$ exactly when
$|a|=1$, we may then assume that $|a| \leq 1.$

We will assume that $|a| <1$ and show that this then implies that $B=
\tri(0,r^*).$

We first note that by Lemma~\REFER{circles}, $B$ is a union
of circles centered at the origin and $B$ contains a sequence of
circles tending to zero.
If $C(0,\rho) \subset B$, then by the forward
invariance of $B$ under $H$, we see that 
$C(0,\rho^j) \subset B.$  
Also we get that if $C(0,\rho) \subset B$,
then by the forward invariance of $B$ under $K$, 
we have $C(0,|a|\rho^m) \subset B.$
Using a change of coordinates $r=\log \rho$
this invariance can be
stated in terms of the new functions
\begin{equation}
t(r)=jr \text{ and } s(r)=mr+c \notag
\end{equation}
where $c=\log |a|<0.$

So the action of $H$ and $K$ on $\tri(0,r^*)$ is replaced
by the action of $t$ and $s$ on $I=[-\infty,\log r^*)$, respectively.  In
particular, we define 
$$B'=\{\log\rho:C(0, \rho)\subset B \} \cup \{ -\infty \}$$  
keeping in mind that $B$ is a union of circles
centered at the origin.  Then 

\begin{align}
& s(B') \subset B'\EQXREF{s(E'')<E''},\\
& s^{-1}(B') \cap I \subset B'\EQXREF{E''backs}, \\
& t(B') \subset B'\EQXREF{t(E'')<E''},\\
& t^{-1}(B') \cap I \subset  B'\EQXREF{E''backt}, \\
& B' \text { is closed in the relative topology on\, } I.
\end{align}

In order to make calculations a bit easier we rewrite 
$s(r)= r_0 + m(r-r_0)$ 
where $r_0=-c/(m-1)>0.$

Hence 
\begin{align}
&s^n(r)=r_0 +m^n(r-r_0), \notag\\ 
&s^{-n}(r)=r_0 +m^{-n}(r-r_0), \notag\\
&t^n(r)=j^nr, \notag\\
&t^{-n}(r) = j^{-n}r.\notag
\end{align}

Consider 
$$(t^{-n}\circ s^{-n}\circ t^n\circ s^n)(r)
=r-r_0+{\frac {r_0} {j^n}}+{\frac {r_0} {m^n}}-{\frac {r_0} {m^n j^n}}.$$

Let 
\begin{equation*}
d_n={\frac {r_0} {j^n}}+{\frac {r_0} {m^n}}-{\frac {r_0} {m^n j^n}}
=r_0 {\frac {m^n+j^n -1} {m^n j^n}}
\end{equation*}
 and note that $0 < d_n \leq r_0$ with $d_n \to 0 $ as $n \to \infty.$ 

We also note that $(t^{-n}\circ s^{-n}\circ t^n\circ s^n)(r)=r-r_0 +d_n$ 
implies that 
$(s^{-n}\circ t^{-n}\circ s^n\circ t^n)(r)=r+r_0 -d_n$ since these two 
functions are inverses of each other.

We claim that $(-\infty, \log r^* -r_0] \subset B'$.  

Let us suppose that this is not the case, and suppose that
$(r',\tilde{r})$ is an interval disjoint from $B'$ with 
$-\infty < r' < \tilde{r} \leq \log r^* -r_0$.  
Since $B'$ is a closed subset of $[-\infty,
\log r^* -r_0]$, we may assume that this interval is expanded so that $r'
\in B'.$  Note that here we used the fact that $B$ contains a
sequence of circles going to 0,
 hence $B'$ contains a sequence of points going to $-\infty.$

Let $r_n'=(t^{-n}\circ s^{-n}\circ t^n\circ s^n)(r')=r'-r_0+d_n.$  
 We claim that each $r_n'$ is in $B'$.  This is almost obvious from
 the invariance of $B'$ under $s$ and $t$ in (\EQREFER{s(E'')<E''}) through
(\EQREFER{E''backt}), but
 some care needs to be taken to insure that each application of $s, \,t,\,
 s^{-1}$, and $t^{-1}$ takes points to the right domain.  
 By (\EQREFER{s(E'')<E''}) we see
 that $s(r'), s^2(r'),\dots, s^n(r') \in B'$.  Hence by (\EQREFER{t(E'')<E''})
we get $(t\circ s^n)(r'), (t^2\circ s^n)(r'), \dots, (t^n\circ s^n)(r')
 \in B'$.  

Since $s^{-1}(r) > r$ for $r \in (-\infty,r_0)$ we see that because 
$(s^{-n}\circ t^n\circ s^n)(r')$ is clearly less than $r'$ (as $t(r)<r$
for $r \in (-\infty,0)$), also each of
$(s^{-1}\circ t^n\circ s^n)(r'), \dots, (s^{-n}\circ t^n\circ s^n)(r')$ must 
be less than $r'<\log r^*$.
Hence by (\EQREFER{E''backs}) we see that each of these points lies in $B'$.

Similarly, since $t^{-1}(r)>r$ for $r \in (-\infty,0)$ and
$(t^{-n}\circ s^{-n}\circ t^n\circ s^n)(r')=r'-r_0+d_n \leq r'<\log
r^* <0$,
also each of
$(t^{-1}\circ s^{-n}\circ t^n\circ s^n)(r'), \dots, 
(t^{-n}\circ s^{-n}\circ t^n\circ s^n)(r')$ 
lies in $I=[-\infty,\log r^*)$.  Hence by (\EQREFER{E''backt}) 
each of these points is in $B'$ and so each $r_n' \in B'$.  

Hence we conclude that $r'-r_0\in B'$ since $B'$ is relatively 
closed in $I$ and $r_n'
\to r'-r_0 \in I.$  Note also that $r_n' \downto r'-r_0.$  

Now we claim that for any $r''\in B'\cap (-\infty, \log r^* -r_0)$, we
have 
$r''+r_0 \in B'$.  Let $r_n''=(s^{-n}\circ t^{-n}\circ s^n\circ t^n)(r'')
=r''+r_0 -d_n$.
Noting that each $r_n'' < r'' + r_0<\log r^*$ we
may again use the invariance of
$B'$ under $s$ and $t$ in (\EQREFER{s(E'')<E''}) through
(\EQREFER{E''backt}) in a similar fashion as above to obtain that
 each $r_n'' \in B'$.   Thus also the limit $r''+r_0 \in B'$.  

Consider again $r_n' \downto
r'-r_0$.  By applying the above claim to each $r_n' \leq r' < \log r^*
-r_0$, we get that each
$r_n'+r_0 \in B'$.  Since $r_n'+r_0 \downto r'$ we then see that
we have contradicted the statement that $(r',\tilde{r})$ is disjoint from
$B'$.

So we conclude that $(-\infty, \log r^*-r_0] \subset B'.$  Clearly then
by the partial backward invariance of $B'$ under the map $t$ we
get
$[-\infty, \log r^*) \subset B'.$  Hence we conclude that $\tri(0,r^*)=B$.
\end{proof}

In order to avoid some technical difficulties we will make use of the
following well known result.

\begin{othertheorem}
A polynomial $f$ of degree $k$ is conjugate near $\infty$
to the map $z \mapsto z^k$ near the origin.  More specifically, there exists
a \nbhd $U$ of $\infty$ such that we have a univalent 
\begin{equation*}
\phi :U \to \tri(0,r^*) \text{ for } 0< r^*<1 \text{ with }
\phi(\infty)=0 \text{ and }
\phi \circ f\circ \phi^{-1}(z) = z^k.
\end{equation*}
\end{othertheorem}

\begin{proof}
After conjugating $f$ by $z \mapsto 1/z$ we may apply Theorem
6.10.1 in ~\cite{Be}, p.~150 to obtain the desired result.
\end{proof}

We will denote the conjugate function of $f$ by $F$, i.e., 
\begin{equation*}
F(z)=\phi\circ f\circ \phi^{-1} (z) = z^k.
\end{equation*}

In order to further simplify some of the following proofs we 
 will assume that $\phi(U)=D=\tri(0,r^*)$.  Note that $U$
 is forward invariant under $f$ since $D=\tri(0,r^*)$ is forward
 invariant under $F$.  We may and 
 will also assume that $U$ is forward invariant
 under $g$ as well.

We now define a corresponding function for $g$ using the same
conjugating map as we did for $f$.  Let  $G$ be the function defined on
$D=\tri(0,r^*)$ given by 
\begin{equation*}
G=\phi\circ g \circ\phi^{-1}.
\end{equation*}
Note that $G(D) \subset D$.  

\vskip.1truein

Via this change of coordinates, we
will use the mappings $F$ and $G$ to obtain information about the 
mappings $f$ and $g$.  In transferring to this simpler coordinate
system we make the following definitions.

Let $W$ be the complement of $E$, i.e., 
\begin{equation*}
W=\CC \setminus E.
\end{equation*}
Note that $W$ is open and it is also completely invariant under both
$f$ and $g$ since it is
the complement of a set which is completely invariant under both $f$
and $g$.

Let $W'$ denote the image of $W$ under $\phi$, i.e., 
$W' = \phi(U \cap W).$  Let $E'$ denote the image of $E$ under $\phi$, i.e.,
$E'=\phi(U \cap E).$
Thus $W'$ is open and $E'$ is closed in the relative
topology of $D$.
Note that $W'$ and $E'$ are disjoint since $W$ and $E$ are disjoint
and $\phi$ is univalent.
Also since $W \cup E = \CC$ it easily follows that $W' \cup E'=
\phi(U) = D$.

By the forward invariance of $W \cap U$ under $f$ we see that
\begin{equation}\label{F(W')<W'}
F(W')=F\circ\phi(W \cap U)=\phi\circ f(W \cap U) \subset \phi (W \cap U)=W'.
\end{equation}

Similarly we get
\begin{equation}\label{F(E')<E'}
F(E') \subset E'.
\end{equation}

  Since $E'$ and $W'$ are disjoint and forward
invariant under $F$, and since $E' \cup W'=D$, we see that 

\begin{align} 
F^{-1}(E') \cap D &\subset E'\EQXREF{E'backF},\\
 F^{-1}(W') \cap D &\subset W'.\EQXREF{W'backF}
\end{align}

Note that in the same
way as we obtained the results for $F$ we get

\begin{align}
G(W') &\subset W'\EQXREF{G(W')<W'},\\
G(E') &\subset E'\EQXREF{G(E')<E'},\\
G^{-1}(E') \cap D &\subset E'\EQXREF{E'backG},\\
G^{-1}(W') \cap D &\subset W'.\EQXREF{W'backG}
\end{align}


\begin{lemma}\label{|a|=1}
If $G(z)=az^l$ with $|a|=1$, then $J_f=J_g$.
\end{lemma}

\begin{proof}
The proof relies on the use of Green's functions.  It is well known
that the unbounded components $F_\infty$ and $G_\infty$ support
Green's functions with pole at $\infty$ which we will denote by $G_f$ and $G_g$
respectively.  It is also well known that on $U$ we have 
\begin{equation*}
G_f(z)=-\log|\phi(z)|
\end{equation*}
since $\phi$ is a map which conjugates $f$ to $z \mapsto z^k$ (see
~\cite{Be}, p.~206).  

Since for $\psi(z)=bz$ where $b^{l-1}=a$, the function $\psi \circ\phi$
conjugates 
$g$ in $U$ to $z \mapsto z^l$, we get in $U$,
\begin{equation*}
G_g(z)=-\log|\psi \circ \phi(z)|=-\log|b\phi(z)|=-\log|\phi(z)|
\end{equation*}
where the last equality uses the fact that $|a|=1$, and so $|b|=1.$

Hence $G_f=G_g$ in $U$.  Since $G_f$ and $G_g$ are each harmonic away
from $\infty$ we get that $G_f=G_g$ on the unbounded component $C$ of 
$F_\infty \cap G_\infty.$  

We claim that this implies that $J_f = J_g$.  Assuming that 
$J_f \neq J_g$, we see by Remark~\REFER{JfneqJg} that there
exists a point which lies in the Julia set of one function, yet
in the unbounded component of the Fatou set of the other function.
Let us therefore suppose that $z_0' \in J_g \cap F_\infty$.  
Let $\gamma$ be a path in
$F_\infty$ connecting $z_0'$ to $\infty.$  We see that $\gamma$
must intersect $\partial C$ somewhere, say at $z_0$.  Since $z_0 \in
F_\infty \cap \partial C$  we get $z_0 \in \partial G_\infty
=J_g$. 

We may
select a sequence $z_n \in C$ such that $z_n \to z_0.$  Since
$z_0$ lies on the boundary of the domain of the Green's function
$G_g$, i.e., $z_0 \in J_g$,
we have $G_g(z_n) \to 0$ (see ~\cite{Be}, p.207).  
Since $z_0$ lies in the domain of the
Green's function $G_f$ we see that $G_f(z_n) \to G_f(z_0)>0.$  We
cannot have both of these happen since $G_f(z_n) = G_g(z_n)$ and so we
conclude that $J_g \cap F_\infty = \emptyset$.  Hence we conclude that
$J_f = J_g$.
\end{proof}


We now are able to prove Theorem~\REFER{main}.

\begin{proof}[Proof of Theorem~\REFER{main}]
Consider whether or not $E$ has nonempty interior.  If $E^o \neq
\emptyset$, then by Lemma~\REFER{Eint} we get $E = \CC$. 

If $E^o
= \emptyset$, then Lemma~\REFER{circles}
implies that the set $W'$ is
a union of open annuli centered
at the origin and hence $E'$ is a union of circles centered at the
origin.  Since $E^o = \emptyset$, the set $E'$ has empty interior.

Since we know by Remark~\REFER{inftynotisol} that there exists a
sequence of points in $E$ tending to infinity when $J_f \neq J_g$, also
$E'$ must contain a corresponding sequence of circles tending to zero.
By Lemma~\REFER{L=azj} we see that the 
function $G$ is of the form
\begin{equation*}
G(z)=az^l
\end{equation*} 
for some non-zero complex number $a$.

By considering the set $E'$, we see that Lemma~\REFER{|a|<1} implies
that $|a|=1.$  We see that Lemma~\REFER{|a|=1} then implies $J_f=J_g$.
\end{proof}

\bibliographystyle{amsplain}
\bibliography{biblio}

\enddocument